\documentclass[12pt]{amsart}

\usepackage{epsfig}
\usepackage{amscd}
\usepackage[mathscr]{eucal}
\usepackage{amssymb}
\usepackage{amsxtra}
\usepackage{amsmath}
\usepackage[all]{xy}

\theoremstyle{plain}

\theoremstyle{definition}

\theoremstyle{remark}

\begin{document}

\bigskip

\title{ On Fermat curves modulo a finite number}

\date{\today}

\author{Yochay Jerby}


%
%
% ...

\begin{abstract} We show that the existence of a non-trivial solution of $x^n+y^n=p^n$, with $p$ a prime number, is equivalent to the existence of a solution of a certain (over-determined) system of $(n-1)$-recursion relations ("zipper" equations) in $\mathbb{Z}_{p-1}$.
\end{abstract}

\maketitle

%
%
\section{Introduction}
\hspace{-0.6cm} The famous Fermat's last theorem states that the equation $x^n + y^n = z^n$ admits no
positive integral solutions, if $n \geq 3$, see \cite{FLT,W} and references therein. But what happens if we ease the condition and
require for points $(x,y) \in \mathbb{Z}^2$ for which $z^n$ just divides $x^n+y^n$, for some fixed $z \in \mathbb{Z}$? In general, this leads us to define the main subject of study in this work:

\bigskip

\hspace{-0.6cm} \bf Definition: \rm For $n \in \mathbb{N}$ and $j \geq 1$ let
$$ T_n(z; j) :=\left \{ (x, y) \vert x^n + y^n \equiv 0 \textrm{ (mod } z^j) \right \} \subset (\mathbb{Z}_{z^j} )^2 $$ be the \emph{$j$-th Fermat tile of radius
$z$}.

\bigskip

\hspace{-0.6cm} We think of elements of $T_n(z ; n)$ as "mock solutions" of the equation $x^n+y^n=z^n$ (as clearly, a genuine solution, if exists, induces an element of $T_n(z ; n)$ but not the other way around). Remarkably, not only do such "mock solutions" exist, they actually satisfy a neat algebraic structure which we describe, using Fermat's little theorem, in section 2. From the elementary features of the plane geometry of Fermat curves, it follows that the tiles $T_n(z ; n)$ satisfy the following property:

\bigskip

\hspace{-0.6cm} \bf Proposition A: \rm If $(x,y,z) \in \mathbb{Z}_+^3$ is a solution of $x^n+y^n=z^n$ then $ (x,y) \in T_n(z ; n) \cap [0,z]^2$.

\bigskip

\hspace{-0.6cm} In section 3 we study the properties of elements $(x,y) \in T_n(p ; n) \cap [0,p]^2$, for $p$ a prime number, and show that such elements are subject to a substantial system of highly non-trivial restrictions. The description of these restrictions requires a study of the functions $$ \begin{array}{ccc} Log_j : \mathbb{Z}^{\ast}_{p^j} \rightarrow \mathbb{Z}_{\phi(p^j) } & ; & Exp_j : \mathbb{Z}_{\phi(p^j) } \rightarrow \mathbb{Z}_{p^j}^{\ast} \end{array}$$ given by $$ \begin{array}{ccc} Log_j(g^s)=[s]_{\phi(p^j)} \in \mathbb{Z}_{\phi(p^j) } & ; & Exp_j(s):=[g^s]_{p^j} \in \mathbb{Z}_{p^j} \end{array}, $$ where $\phi(p^j)=(p-1) \cdot p^{j-1}$ is the totinet function and $j \geq 1$. The main feature is the definition of a natural class of functions of the form $A^a_j : \mathbb{Z}_{p-1} \times \mathbb{Z}_{p^{j-1}} \rightarrow \mathbb{Z}_p$, for $ a \in \mathbb{Z}_{p-1}$. In terms of these functions, the restrictions are given as follows:

\bigskip

\hspace{-0.6cm} \bf Theorem B \rm ("zipper" relations): If $(x,y) \in T_n(p ; n) \cap [0,p]^2$ there exists an element $s = s_{(x,y)} \in \mathbb{Z}_{p-1}$ and $a \in \mathbb{Z}_{p-1}$ which satisfy the following overdetermined double recursion (\emph{"zipper"}) relations $$ \left \{ \begin{array}{c} r_0=A^0_1(s)= A^a_1(s) \\ r_1 = A^0_2(s ; r_0) = A^a_2(s ; r_0) \\ \vdots \\ r_{n-1} = A_n^0(s ; r_0 ...,r_{n-2}) = A_n^a(s ; r_0 ...,r_{n-2}) \end{array} \right. $$

\bigskip

\hspace{-0.6cm} As one can see the number of equations in the variable $s \in \mathbb{Z}_{p-1}$ grows with $n$ (in the Pythagorean case $n=2$, there is one equation). In particular, Fermat's last theorem would follow from showing that the zipper relations have no solutions for $n \geq 3$. We define and study various properties of the zipper relations in section 4.

\bigskip

\hspace{-0.6cm} The rest of the work is organized as follows: In section 2 we describe Fermat tiles, in section 3 we study $Log_j$ and $Exp_j$ and define $A^a_j$. In section 4 we define the zipper relations.

\section{The geometric structure of Fermat tiles}

\hspace{-0.6cm} Before describing Fermat tiles in general, let us start with a few examples (which justify the term \emph{tile}).

\bigskip

\hspace{-0.6cm} \bf Example 2.1 \rm (n = 2): Figure 4 shows the first Fermat tile $T_2(5, 1)$ of radius $5$.

\begin{figure}[!ht]
\includegraphics[scale=0.4]{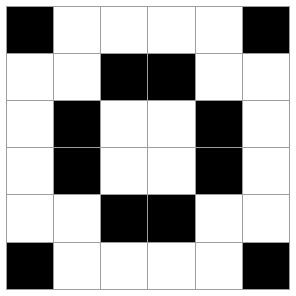}
\caption{ First Fermat tile $T_2(5, 1)$ of radius $5$. \label{overflow}}
\end{figure}

\hspace{-0.6cm} Figure 5 shows the second Fermat tile $T_2(5, 2)$ of radius $5$.
\begin{figure}[!ht]
\includegraphics[scale=0.4]{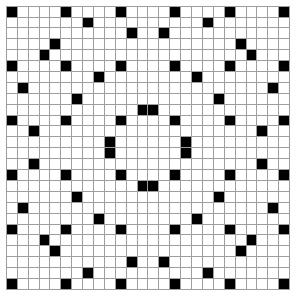}
\caption{Second Fermat tile $T_2(5, 2)$ of radius $5$. \label{overflow}}
\end{figure}

\hspace{-0.6cm} Note that, as expected, one has $(3, 4), (4, 3),(5,0),(0,5) \in T_2(5, 2) \cap [0, 5]^2$.

\bigskip

\hspace{-0.6cm} \bf Example 2.2 \rm ($n \geq 3$): Figure 6 shows the first Fermat tile $T_4(7, 1)$ of radius $7$:
\begin{figure}[!ht]
\includegraphics[scale=0.4]{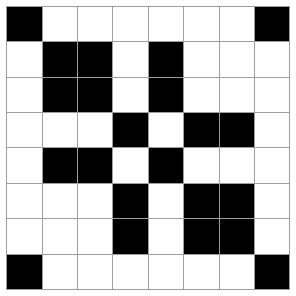}
\caption{Second Fermat tile $T_3(7, 1)$ of radius $7$. \label{overflow}}
\end{figure}

\hspace{-0.6cm} Figure 7 shows the first Fermat tile $T_4(17, 1)$ of radius $17$

\begin{figure}[!ht]
\includegraphics[scale=0.35]{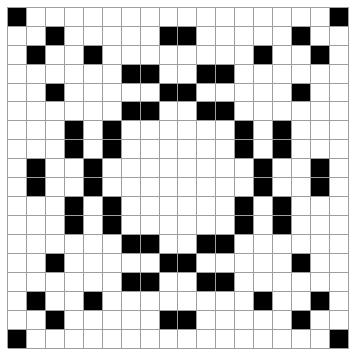}
\caption{First Fermat tile $T_4(17, 1)$ of radius $17$. \label{overflow}}
\end{figure}

\hspace{-0.6cm} We refer to
$$C_n(z) := \left \{ (x,y) \vert x^n + y^n = z^n \right \} \subset \mathbb{R}^2 $$
as the $n$-th Fermat curve of radius $ z \in \mathbb{R}$. Note that Fermat's last theorem is equivalent to stating that $C_n(z) \cap \mathbb{N}^2 = \left \{ (z,0) , (0,z) \right \}$ for any $n \geq 3$ and $z \in \mathbb{N}$. Let us proceed with the following remark:

\bigskip

\hspace{-0.6cm} \bf Remark 2.3 \rm (plane geometry of Fermat curves): Figure 5 shows $C_2(5)$, the circle of radius $z = 5$, with
the integer lattice:

\begin{figure}[!ht]
\includegraphics[scale=0.27]{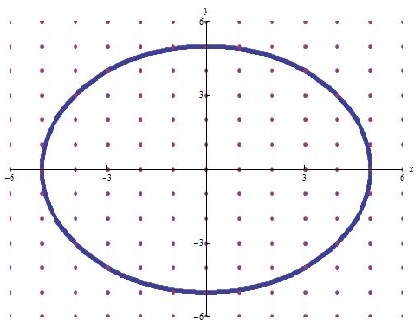}
\caption{Graph of $C_2(5)$ together with the integer lattice $\mathbb{Z}^2$. \label{overflow}}
\end{figure}

\hspace{-0.6cm} Recall that a solution $(x, y, z) \in \mathbb{N}^3$ of $x^2 + y^2 = z^2$ is
called a Pythagorean triple. Note that $C_2(5) \cap \mathbb{N}^2 = \left \{ ( 3, 4), ( 4, 3) , ( 5 ,0 ) , (0, 5) \right \}$.
Figure 6 shows $C_8(5)$ with the integral lattice

\begin{figure}[!ht]
\includegraphics[scale=0.37]{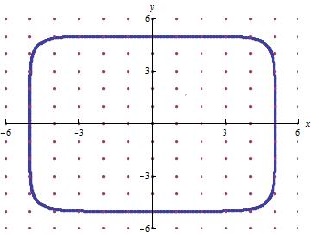}
\caption{Graph of $C_8(5)$ together with the integer lattice $\mathbb{Z}^2$. \label{overflow}}
\end{figure}

\hspace{-0.6cm} Figure 7 shows $C_9(10)$ with the integer lattice

\begin{figure}[!ht]
\includegraphics[scale=0.27]{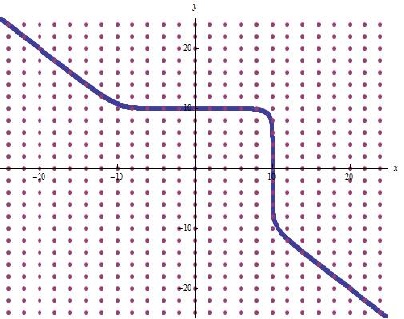}
\caption{Graph of $C_9(10)$ together with the integer lattice $\mathbb{Z}^2$. \label{overflow}}
\end{figure}

\hspace{-0.6cm} In particular, it is easy to see that Fermat curves satisfy $C_n(z) \cap (\mathbb{R}^+)^2 \subset [0,z]^2$ for any $n \in \mathbb{N}$ and $z \in \mathbb{R}^+$.

\bigskip
\hspace{-0.6cm} Note that, in view the above remark, showing that $$ T_n(z ;n) \cap [0,z]^2 = \left \{ (0,0),(z,0),(0,z),(z,z) \right \} $$ would imply Fermat's last
theorem. Let us now turn to describe the structure of the $j$-th Fermat tile, $T_n(p, j)$, for $p$ a prime number. Let
us start with the following:

\bigskip

\hspace{-0.6cm} \bf Lemma 2.4: \rm Let $S_n(p) \subset \mathbb{ Z}_p$ be the solution set of the equation $x^n + 1 \equiv 0 \textrm{ (mod } p)$ and let
$g \in \mathbb{Z}_p$ be a primitive generator.

\bigskip

(a) If $2n \mid (p- 1)$ then $S_n(p) = \left \{g^{a_0+i m} \right \}_{i=0}^{n-1}$ where $a_0 := \frac{p-1}{2n}$
and $m := \frac{p-1}{n}$.

\bigskip

(b) If $2n \nmid (p - 1)$ then $S_n(p) = \emptyset$.

\bigskip

\hspace{-0.6cm} \bf Proof: \rm As $ g \in Z_p$ is a primitive generator we have
$$ \left \{ g^s \vert s=0,...,p-2 \right \} \subset \mathbb{Z}^{\ast}_p$$
By Fermat's little theorem $g^{ \frac{p-1}{2}} \equiv -1 \textrm{ (mod }p)$. Hence, we are looking for some $x = g^a$
such that $$ g^{n \cdot a} \equiv g^{\frac{p-1}{2}} \textrm{ (mod }p).$$
Again, by Fermat's little theorem, this is equivalent to solving $$ n \cdot a \equiv \frac{p - 1}{2} \textrm{ (mod } p- 1).$$
The solutions are given by $a_i := a_0+i m$ for $ i = 0,..,n-1$. $\hspace{5cm} \square$

\bigskip

\hspace{-0.6cm} Consider the following example:

\bigskip

\hspace{-0.6cm} \bf Example 2.5: \rm Let $n = 4$ and $p = 17$. One has $a_0 = 2,m = 4$ and $g = 3$. Hence
$$S_4(17) = \left \{3^2, 3^6, 3^{10}, 3^{14} \right \} = \left \{9, 15, 8, 2 \right \} \subset \mathbb{ Z}_{17},$$
which coincides with the first row of Fig. 4, as expected.

\bigskip

\hspace{-0.6cm} We have:

\bigskip

\hspace{-0.6cm} \bf Proposition 2.6 \rm (reducibility): For any $1 \leq j \leq n$ the following holds:

\bigskip

(a) If $2n \mid (p - 1)$ then $T_n(p, j) = \left \{ \left ( a,a x^{p^{j-1}} \right ) \vert x \in S_n(p) , a \in \mathbb{Z}_{p^j} \right \} \cup \left ( p \mathbb{Z}_{p^{j-1}} \right )^2$.

\bigskip

(b) If $2n \nmid (p - 1)$ then $T_n(p, j) = \left ( p \mathbb{Z}_{p^{j-1}} \right )^2$.
\bigskip

\hspace{-0.6cm} For $x \in S_n(p), j \in \mathbb{N}$ consider the linear function $f^j_x: \mathbb{Z}_{p^j} \rightarrow \mathbb{Z}_{p^j}$ given by $a \mapsto a x^{p^{j-1}}$. Set $$ T(x;j) := graph(f^j_x)=\left \{ \left (a, a x^{p^{j-1}} \right ) \vert a \in \mathbb{Z}_{p^j} \right \}.$$ Consider the following example:

\bigskip

\hspace{-0.6cm} \bf Example 2.7: \rm Let $n=4$ and $p=17$. Figure 8 shows the graphs of $f^1_x : \mathbb{Z}_{17} \rightarrow \mathbb{Z}_{17}$ for $x= 9,15,8,2$

\begin{figure}[!ht]
\includegraphics[scale=0.31]{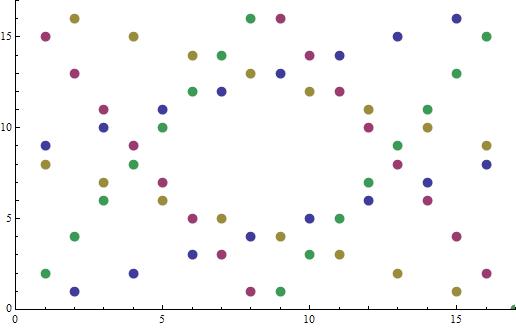}
\caption{Graphs of $f^1_x$ for $x = 9,15,8,2$. \label{overflow}}
\end{figure}

\hspace{-0.6cm} As one can see, Figure 8 coincides with the interior of the Fermat tile $T_4(17,1)$, presented in Figure 7. Figure 9 presents the four linear components separately:

\begin{figure}[!ht]
\includegraphics[scale=0.3]{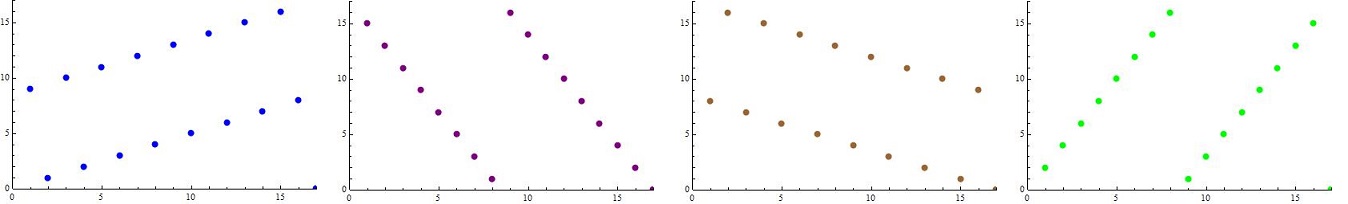}
\caption{Separated graphs of $f^1_x$ for $x = 9,15,8,2$. \label{overflow}}
\end{figure}

\hspace{-0.6cm} Figure 10 shows the graphs of $f^2_x : \mathbb{Z}_{17^2} \rightarrow \mathbb{Z}_{17^2}$ for $x= 9,15,8,2$
\begin{figure}[!ht]
\includegraphics[scale=0.31]{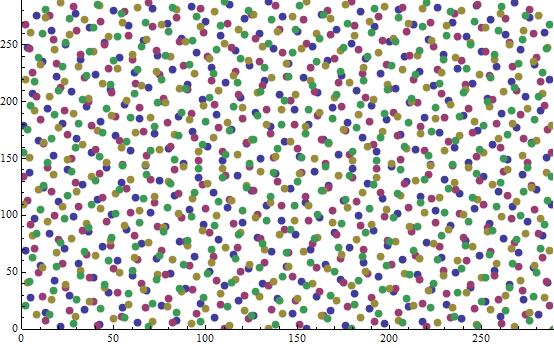}
\caption{Graphs of $f^2_x$ for $x = 9,15,8,2$. \label{overflow}}
\end{figure}

\newpage

\hspace{-0.6cm} Figure 11 presents the four linear components separately:

\begin{figure}[!ht]
\includegraphics[scale=0.325]{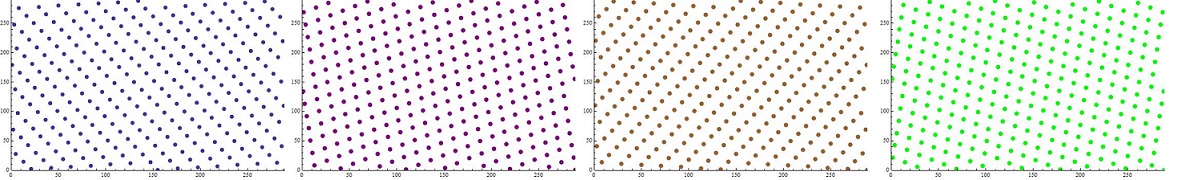}
\caption{Separated graphs of $f^2_x$ for $x = 9,15,8,2$. \label{overflow}}
\end{figure}

\hspace{-0.6cm} In view of the above, Fermat's last theorem would thus follow from showing $$ T(x ; n) \cap [0,p]^2 = \left \{ (0,0) \right \} \hspace{0.5cm} \textrm{ for all } \hspace{0.25cm} x \in S_n(p).$$ That is, showing $f_x^n(a)>p$ for any $1 \leq a \leq p$ . First, let us make the following remark:

\bigskip

\hspace{-0.6cm} \bf Remark 2.8 \rm (empirics): $T(x ;j)$, is the graph of a line of slope $x^{p^{j-1}}$ in $\mathbb{Z}^2_{p^j}$, for $x \in S_n(p)$. In parctice, such a line actually traces a $2$-dimensional lattice in $\mathbb{Z}_{p^j}^2$, due to the truncation caused by the quotient relation (see Fig. 11 and 13). Empirics show that, for various values of $j$ and $x$, typically, the values of $T(x;j)$ can be bounded by a line $ a_2 =\frac{m_j(x) \cdot a_1}{n_j(x)} + \frac{p^j}{n_j(x)}$ with $m_j(x)<\sqrt{p}^j$. Figure 12 shows $(a, f^3_x(a))$ for $9 \in S_4(17)$ in $[0,p^2] \times [0,p^3]$ and $[0,p^2] \times [0,p^2]$, together with the bounding line $a_2 = -\frac{23}{47} \cdot a_1 + \frac{17^3}{47}$:

\begin{figure}[!ht]
\includegraphics[scale=0.3]{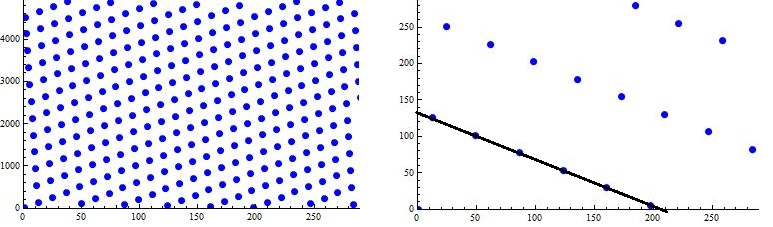}
\caption{$(a,f^3_x(a))$ in $[0,17^2] \times [0,17^3]$ and $[0,17^2] \times [0,17^2]$. \label{overflow}}
\end{figure}
\hspace{-0.6cm} Note that, if such a line exists in general, for $j \geq 3$, it would need to go more than $p$ steps in the $a_1$-axis to go below $p$ in the $a_2$-axis and, in particular, $f^j_x(a)>p$ if $0 \leq a \leq p$. However, showing that such a line exists, in general, requires for a more solid understanding of $T(x;j) \subset \mathbb{Z}_{p^j}^2$ (the subject of the following section).

\hspace{-0.6cm} It is also interesting to note the following: set $x_i = g^{a_i} \in S_n(p)$ for $i=0,...,n-1$ and denote $\theta_i^j(p,n)$ be the number of elements $0 \leq a \leq p $ such that $f^j_{x_i}(a)<p$ (we want $\theta^n_i(p,n)=0$ for all $i$). For $n=4$ we actually have $$ \tiny \left ( \begin{array}{ccccccccccccccccccccc}
p & 17 & 41 & 73 & 89 & 97 & 113 & 137 & 193 & 233 & 241 & 257 & 281 & 313 & 337 & 353 & 401 \\
\theta^2_1 & 1& 3& 2& 0& 3& 0& 1& 2& 1& 0& 1& 1& 1& 1& 1& 0\\
\theta^2_2 & 1& 0& 0& 1& 0& 1& 0& 0& 0& 1& 1& 0& 1& 1& 1& 3\\
\theta^2_3 & 1& 0& 0& 1& 0& 1& 0& 0& 0& 1& 1& 0& 1& 1& 1& 3\\
\theta^2_4 & 1& 3& 2& 0& 3& 0& 1& 2& 1& 0& 1& 1& 1& 1& 1& 0 \\
\sum \theta_i^3 & 0 & 0 & 0 & 0& 0 & 0& 0 & 0 & 0 & 0& 0 & 0 & 0 & 0 & 0 & 0 \end{array} \right ) .
$$

\section{Some remarks on the arithemetics of the ring of invertibles $\mathbb{Z}^{\ast}_{p^j}$}

\hspace{-0.6cm} Let $g \in \mathbb{Z}_p$ be a fixed primitive generator. It is easy to see that $g$ is a primitive generator for $\mathbb{Z}_{p^j}$ for all $j \geq 1$, as well. Recall that the ring of invertibles is given by $$\mathbb{Z}^{\ast}_{p^j} = \left \{ a \vert a \not \equiv 0 \textrm{ (mod } p) \right \} \subset \mathbb{Z}_{p^j}.$$ Set $\phi(p^j):= (p-1) \cdot p^{j-1}$. Let us consider the following two functions $$ \begin{array}{ccc} Log_j : \mathbb{Z}^{\ast}_{p^j} \rightarrow \mathbb{Z}_{\phi(p^j) } & ; & Exp_j : \mathbb{Z}_{\phi(p^j) } \rightarrow \mathbb{Z}_{p^j}^{\ast} \end{array}$$ given by $Log_j(g^s)=[s]_{\phi(p^j)} \in \mathbb{Z}_{\phi(p^j) }$ and $Exp_j(s):=[g^s]_{p^j} \in \mathbb{Z}_{p^j}$. By definition, the amobea corresponding to the linear tile $T(x ; j)$ is simply given by the following affine line $$ \mathcal{A}(x_i ; j):= Log_j \left ( T(x_i ; j) \cap (\mathbb{Z}_{p^j}^{\ast})^2 \right ) = \left \{ (s_1,s_2) \vert s_1-s_2=a_i \cdot p^{j-1} \right \} \subset \mathbb{Z}^2_{\varphi(p^j)}, $$ where $Log_1(x_i) = a_i = \frac{ (1+2i) \cdot (p-1) }{2n} $, with $i=0,...,n-1$. In particular, we want to show $$max \left\{ Exp_n(s), Exp_n(s+a_i \cdot p^{n-1}) \right \} \geq p,$$ for any $s \in \mathbb{Z}_{\phi(p^n)}$ and $i=0,...,n-1$.

\bigskip

\hspace{-0.6cm} \bf Remark 3.1 \rm (Transcendentality of $Log_j$ and $Exp_j$): It should be noted that finding a general \emph{full} explicit description of $Log_j$ and $Exp_j$ is considered a completely transcendental question. For instance, Fig. 13 is a graph of $Log_1(a) \in \mathbb{Z}_{96}$ for $a \in \mathbb{Z}_{97}$ with $g=5$:

\begin{figure}[!ht]
\includegraphics[scale=0.4]{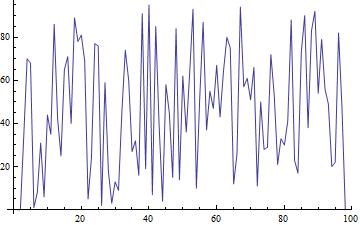}
\caption{$Log_1(a)$ for $p=97$ and $g=5$. \label{overflow}}
\end{figure}

\bigskip

\hspace{-0.6cm} Note that, by definition, $[Exp_{j+1}(s)]_{p^j}=Exp_j([s]_{\phi(p^j)})$. Hence, for any $j \geq 1$ we can define the function $g_j : \mathbb{Z}_{\varphi(p^{j+1})} \rightarrow \mathbb{Z}_p$ given by $$ g_j(s) : = \frac{ Exp_{j+1}(s) - Exp_j ( [s]_{\phi(p^j)}) } {p^j}.$$ Consider the formal series $$ Exp(s,T):=Exp_1([s]_{p-1}) + \sum_{i=1}^{\infty} g_j([s]_{\phi(p^{j+1})}) \cdot T^j.$$ The formal series satisfies $Exp_j(s) = [Exp(s,p)]_{p^j}$. In particular, note that $$ \left \{ \begin{array}{c} Exp_j(s) \leq p \Leftrightarrow g_1(s) =...=g_{j-1}(s)=0, \\ Exp_j(s+a \cdot p^{j-1} ) \leq p \Leftrightarrow g_1(s+a \cdot p) =...=g_{j-1}(s+ a \cdot p^{j-1} )=0 \end{array} \right . $$ Note that an element $r \in \mathbb{Z}_{p^j}$ can be uniquely expressed as $$ r =r_0 + r_1 \cdot p + ... + r_{j-1} \cdot p^{j-1},$$ such that $r_i \in \mathbb{Z}_p$ for any $1 \leq i \leq j-1$. For $ a \in \mathbb{Z}_{p-1}$, let $h^a_j : \mathbb{Z}_{p-1} \times \mathbb{Z}_{p^{j}} \rightarrow \mathbb{Z}_{p}$ be the function given by $$ h_j^a(s ; r_0,...,r_{j-1}) := g_j \left (s ,r_0 + r_1 \cdot p +... + r_{j-1} \cdot p^{j-1} + a \cdot p^j \right ).$$ It is easy to see that:

\bigskip

\hspace{-0.6cm} \bf Lemma 3.2 \rm (shift): For any $a \in \mathbb{Z}_{p-1}$ and $j \geq 1$ the following holds $$ h_j^a(s ; r_0,...,r_{j-1}) = h^0_j(s+a ; r_0 + a+r(s+a) ,...,r_{j-1} + a+r(s+a)), $$ where $r(m)=\frac{(m-[m]_{p-1})}{p-1}$.

\bigskip

\hspace{-0.6cm} For instance, consider the following example:

\bigskip

\hspace{-0.6cm} \bf Example 3.3 \rm ($p=5$ and $j=1$): The following is a table of the values of $h^0_1(s, r)$ for $p=5$
$$ \left ( h^0_1 (s,r) \right )_{0,0}^{3,4} = \left ( \begin{array}{ccccc}
0 & 3 & 1 & 4 & 2 \\
0 & 1 & 2 & 3 & 4 \\
0 & 2 & 4 & 1 & 3 \\
1 & 0 & 4 & 3 & 2 \end{array} \right ).$$
The values of $h^3_1(s_1,s_2)$ are given by $$\left ( h^3_1 (s,r) \right )_{0,0}^{3,4}= \left ( \begin{array}{cccccc}
3 & 2 & 1 & 0 & 4 \\
2 & 0 & 3 & 1 & 4 \\
4 & 0 & 1 & 2 & 3 \\
3 & 0 & 2 & 4 & 1 \end{array} \right )$$ Note that $(s^0,r^0)=(3,1)$ is a point satisfying $h_1^0(s^0,r^0)=h_1^3(s^0,r^0)=0$. Indeed, $$ \begin{array}{ccc} Exp_2(7)=Exp_2(3+1 \cdot 4) = 3 & ; & Exp_2(2) = Exp_2([22]_{20})=Exp_2(7+ 3 \cdot 5) = 4 \end{array}, $$ which represents the Pythagorean point $(3,4) \in C_2(5)$, as expected.

\bigskip

\hspace{-0.6cm} In general, for any $j \geq 1$, set $$Z_j(a) : = \left \{ (s,r) \vert h_k^a(s,r_0,...,r_{k-1}) = 0 \textrm{ for any } k \leq j \right \} \subset \mathbb{Z}_{p-1} \times \mathbb{Z}_{p^{j}}. $$ Let $p$ be a prime and $2 \leq n \in \mathbb{N}$ an integer such that $2n /(p-1)$. Set $ a_i := \frac{ (1+2i)(p-1)}{2n}$ for $i=1,...,n-1$. Our main object of study now is $ Z_{n-1}(0) \cap Z_{n-1}(a_i)$. In particular, note that according to the above, showing $$ Z_{n-1}(0) \cap Z_{n-1}(a_i) = \emptyset $$ would imply Fermat's last theorem. Consider, for instance, the following example:

\bigskip

\hspace{-0.6cm} \bf Example 3.4 \rm ($p=5$ and $j=2$): Note from the previous example that, for $(s,r_0)=(3,1)$, we have
$$ h^0_1(3 ; 1)=h^3_1(3;1) = h^0_1(6 ; 4) =h^0_1(1,0)=0.$$ Note also that $$ \begin{array}{ccc} h_1^0(3,r) =[4 \cdot (r-1) ]_5 & ; & h_1^0(2,r)=[2 \cdot r]_5 \end{array}. $$ For $j=2$, we have
$$ \left ( h_2^0 \left (3 ; r_0,r_1 \right ) \right )_{0,0}^{4,4}= \left ( \begin{array}{ccccc}
0 & 4 & 3 & 2 & 1\\
0 & 4 & 3 & 2 & 1\\
1 & 0 & 4 & 3 & 2\\
0 & 4 & 3 & 2 & 1\\
1 & 0 & 4 & 3 & 2 \end{array} \right ), $$ and $$ \left ( h_2^0 \left (2 ; r_0,r_1 \right ) \right )_{0,0}^{4,4}= \left ( \begin{array}{ccccc}
0 & 2 & 4 & 1 & 3\\
2 & 4 & 1 & 3 & 0\\
0 & 2 & 4 & 1 & 3\\
0 & 2 & 4 & 1 & 3\\
0 & 2 & 4 & 1 & 3 \end{array} \right ). $$ As one can see, we can still express $$ \begin{array}{ccc} h_2^0(3,r_0,r_1) =[4 \cdot (r_1-A^0_2(3; r_0)) ]_5 & ; & h_2^0(2,r_0,r_1)=[2 \cdot (r_1-A^0_2(2 ;r_0)) ]_5 \end{array}, $$ where $A^0_2: \mathbb{Z}_{p-1} \times \mathbb{Z}_p \rightarrow \mathbb{Z}_p$ is a non-linear function, determining the position of the zero in the $r_0$-row (see Lemma 3.5 below). In particular, note that the linear system $$\begin{array}{ccc} h^0_2 (3 ; 1,r_1) =[ 4 \cdot r_1 ]_5=0 & ; & h^0_2(2 ; 1 ,r_1)=[2 \cdot (r_1-4) ]_5= 0 \end{array}$$ has no solution. Hence $Z_1(0) \cap Z_1(3) = \left \{ (3,1) \right \}$ while $Z_2(0) \cap Z_2(3) = \emptyset$. Finally, it should be noted that the values of $Exp_1(s)$ for $s \in \mathbb{Z}_4$ are given by $$ \left ( \begin{array}{ ccccc} s & 0 & 1 & 2 & 3 \\ Exp_1(s) & 1 & 2 & 4 & 3 \end{array} \right ). $$ In particular, $h_2^0$ can be expressed in the following form $$h_2^0(s;r_0,r_1) =[Exp_1([s+3]_4) \cdot (r_1-A^0_2(s ; r_0)) ]_5 . $$

\bigskip

\hspace{-0.6cm} In general, we have:

\bigskip

\hspace{-0.6cm} \bf Lemma 3.5: \rm For any $ a \in \mathbb{Z}_{p-1}$ and $j \geq 1$ there exists $A^a_j : \mathbb{Z}_{p-1} \times \mathbb{Z}_{p^{j-1}} \rightarrow \mathbb{Z}_p$ and $s_0 \in \mathbb{Z}_{p-1}$, such that $$ h^a_j (s ; r_0,...,r_{j-1}) = [Exp_1(s+a+s_0) \cdot ( r_{j-1}- A_j^a(s; r_0,..,r_{j-2}))]_p. $$

\bigskip

\hspace{-0.6cm} First, note that, due to the shift property $$ A_j^a(s ; r_0 , ...,r_{j-2}) = A_j^0(s+a ; r_0 + a + r(s+a) ; ... ; r_{j-2}+a + r(s+a) ) - (a+r(s+a)). $$ For instance, consider the following example:
\bigskip

\hspace{-0.6cm} \bf Example 3.6: \rm For $p=5$ and $j=1$ we have $$ \left ( \begin{array}{ccccc} s & 0 & 1 & 2 & 3 \\
A_1^0 & 0 & 0 & 0 &1 \\ A_1^3 & 3 & 1 & 1 & 1 \end{array} \right ). $$ Indeed, note that $$\begin{array}{ccc} A^3_1(0) = [A_1^0(3) -3]_5= [1-3]_5 =[-2]=3 & ; & A^3_1(1) = [A^0_1(4)-(3+r(1+3))]_5= [-4]_5=1 \end{array}. $$

\hspace{-0.6cm} The study of various properties of the functions $A^a_j$, defined in Lemma 3.5, is the subject of the next section.

\section{The double recursion ("zipper") equations}

\hspace{-0.6cm} For any $1 \leq k \leq j$, set $$Z^k_j(a) : = \left \{ (s,r) \vert h_k^a(s,r_0,...,r_{k-1}) = 0 \right \} \subset \mathbb{Z}_{p-1} \times \mathbb{Z}_{p^{j}}. $$ Note that $$Z_j(a) = Z_j^1(a) \cap ... \cap Z_j^j(a).$$ By definition, $$ Z_j^k(a) = \left \{ (s,r_0,...,r_{k-2}, A_k^a(s ; r_0, ...,r_{k-2}),r_k,...,r_{j-1} )\right \} \subset \mathbb{Z}_{p-1} \times \mathbb{Z}_{p^j}. $$ Hence, from the above description, we deduce that $(s ; r_0, ...,r_{j-1} ) \in Z_j(a)$ if and only if it is a solution of the following recursion relations $$ \begin{array}{ccccccc} r_0= A^a_1(s) & ; & r_1 = A^a_2(s ; r_0) & ; & ...& ; & r_{j-1} = A_j^a(s ; r_0 ...,r_{j-2}) \end{array} $$ In conclusion, combining the above relations for $Z_j(0)$ and $Z_j(a)$, we get:

\bigskip

\hspace{-0.6cm} \bf Theorem 4.1 \rm ("zipper" relations): Let $j \geq 1$. Then $(s ; r_0, ...,r_{j-1} ) \in Z_j(0) \cap Z_j(a)$ if and only if it is a solution of the following (overdetermined) double recursion relations $$ \left \{ \begin{array}{c} r_0=A^0_1(s)= A^a_1(s) \\ r_1 = A^0_2(s ; r_0) = A^a_2(s ; r_0) \\ \vdots \\ r_{j} = A_{j+1}^0(s ; r_0 ...,r_{j-1}) = A_{j+1}^a(s ; r_0 ...,r_{j-1}) \end{array} \right. $$

\bigskip

\hspace{-0.6cm} As mentioned in section 3, showing that the zipper relations admit no solutions for $j=n-1$ would imply Fermat's last theorem. Let us conclude this section with a few further remarks on the properties of the zipper relations.
\bigskip

\hspace{-0.6cm} \bf Remark 4.2 \rm (the geometry of $A^a_j$): The zipper relations are given in terms of the functions $A^a_j(s ;r_0,...,r_{j-2} )$ which, by definition, represent a parametrization of the zero set of the function $h^a_j(s ; r_0,..., r_{j-1})$. It is interesting to note that the dependency of $A^a_j$ on the $s$-parameter is essentially different than its dependency on the $r=(r_0,...,r_{j-2})$ parameters. For instance, let $p=97$, $n=4$ and $a= (96/8)=12$. Figure 14 shows a matrix-plot of the function $h^a_2( s_0 ; r_0 ,r_1)$ for $s_0=Log_1(11)$ fixed.

\begin{figure}[!ht]
\includegraphics[scale=0.4]{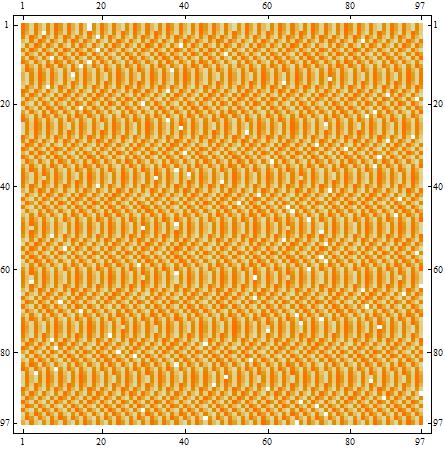}
\caption{Matrix-plot for $h^a_2(s ; r_0,r_1)$ for $(r_0,r_1) \in \mathbb{Z}_p^2$. \label{overflow}}
\end{figure}

\hspace{-0.6cm} In particular, it might seem that the zeros of $h_j^0(s ;r_0,r_1)$ (white points) are randomly distributed in the $r$-coordinates. However, they are determined as the local minima of a one-parametric family of quadrics dominating the picture. Hopefully, we would give a more detailed description of this one parametric family of quadrics in a future work.

\hspace{-0.6cm} On the other hand, it seems that, in the $s$-coordinate, the function $A^a_j(s ; r)$ inherits the chaotic behavior of $Log_1(s)$. For instance, Figure 15 is a graph of $A_1^a(s)$ for the same parameters as above:

\begin{figure}[!ht]
\includegraphics[scale=0.45]{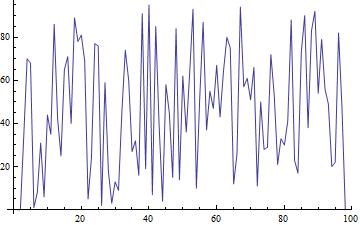}
\caption{$A^a_1(s)$ for $s \in \mathbb{Z}_{96}$ with $a=12$. \label{overflow}}
\end{figure}

\hspace{-0.6cm} Furthermore, it is also interesting to note, for the above parameters, that the solutions of the first zipper equation are give by $$ \left \{ s \vert A^0_1(s) = A^a_1(s) \right \} =\left \{ Log_1(11),Log_1(22),Log_1(33) \right \} \subset \mathbb{Z}_{96}. $$ On the other hand, the second zipper equations for these elements are already non-zero $$ \left ( A_2^0(Log_1(11) , A^0_1(Log_1(11))-A_2^a(Log_1(11) , A^0_1(Log_1(11)) \right ) =$$
$$= \left ( A_2^0(Log_1(22) , A^0_1(Log_1(22))-A_2^a(Log_1(22) , A^0_1(Log_1(22)) \right ) =$$
$$= \left ( A_2^0(Log_1(33) , A^0_1(Log_1(33))-A_2^a(Log_1(33) , A^0_1(Log_1(33)) \right ) = 3 \neq 0.$$ Moreover, note that the value of the second zipper equation is actually independent of the element of the solution set of $A^a_1 (s)= A^0_1(s)$ chosen. This repeats itself in other cases as well.

\bigskip

\hspace{-0.6cm} The above remark shows that $A^a_j(s,r)$ seems to be chaotic in the $s$-parameter and non-linear in the $r$-parameter. In order to overcome this, let us note that if $(s ; A^0_1(s) , A^0_2(s ; A_1^0(s)))$ is a solution, for instance, of the first two zipper equations it needs to satisfy $$ \begin{array}{ccc} h_1^a \left (s ; A^0_1(s) \right ) = 0 & ; & h_1^0 \left (s ; A^a_1(s) \right ) = 0 \\
h_2^a \left (s ; A^0_1(s), A^0_2(s ; A_1^0(s) ) \right )=0 & ; & h_2^0 \left (s ; A^a_1(s), A^a_2(s ; A_1^a(s) ) \right )=0 \end{array}. $$ In view of this, let us define for $i=1,2$ the functions $H^0_i,H^a_i : \mathbb{Z}_{p-1} \rightarrow \mathbb{Z}_p$ given by $$\begin{array}{ccc}
H_1^a(s):= h_1^a \left (s ; A^0_1(s) \right ) & ; & H_1^0(s):= h_1^0 \left (s ; A^a_1(s) \right ) \\
H^a_2(s) :=h_2^a \left (s ; A^0_1(s), A^0_2(s ; A_1^0(s)) \right ) & ; & H^0_2(s) :=h_2^0 \left (s ; A^a_1(s), A^a_2(s ; A_1^a(s)) \right )
\end{array} . $$ Further define $\widetilde{H}^0_i,\widetilde{H}^a_i: \mathbb{Z}_{p}^{\ast} \rightarrow \mathbb{Z}_p$ by $\widetilde{H}^0_i (t) = H^0_i(Log_1(t))$ and $\widetilde{H}^a_i (t) = H^a_i(Log_1(t))$. Remarkably, contrary to $A^0_i(s ; r)$ which are chaotic as functions of $s$, the functions $\widetilde{H}^0_i,\widetilde{H}^a_i$ are actually quite patterned. In fact, Figure 16 shows the graphs of $\widetilde{H}^0_1,\widetilde{H}^a_1$ for $p=97$ and $a=12$:

\begin{figure}[!ht]
\includegraphics[scale=0.45]{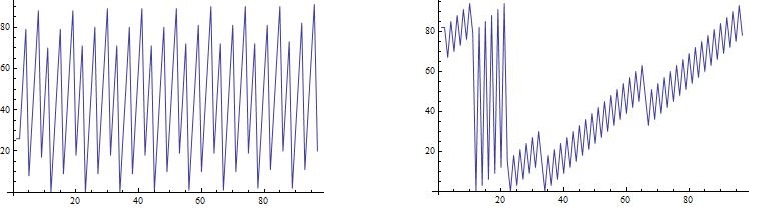}
\caption{$\widetilde{H}^0_1(t)$ (left) and $\widetilde{H}^a_1(t)$ (right) for $t \in \mathbb{Z}^{\ast}_{97}$ with $a=12$. \label{overflow}}
\end{figure}

\hspace{-0.6cm} Note that, as expected, $ \widetilde{H}^0_1(t)=\widetilde{H}^a_1(t)=0$ for $t=11,22,33$. Figure 17 shows the graphs of $\widetilde{H}^0_2(t),\widetilde{H}^a_2(t)$:

\begin{figure}[!ht]
\includegraphics[scale=0.45]{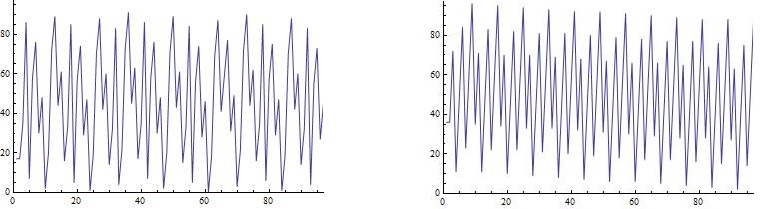}
\caption{$\widetilde{H}^0_2(t)$ (left) and $\widetilde{H}^a_2(t)$ (right) for $t \in \mathbb{Z}^{\ast}_{97}$ with $a=12$. \label{overflow}}
\end{figure}

\hspace{-0.6cm} As expected, one clearly sees that these four graphs have no common zero. In order to further describe $\widetilde{H}^a_i$ let us introduce their \emph{derivatives}. Let $F: \mathbb{Z}_p^{\ast} \rightarrow \mathbb{Z}_p$ be
a function. We refer to the function $D(F): \mathbb{Z}_p^{\ast} \rightarrow \mathbb{Z}_p$, given by $D(F)(t) :=[f(t+1)-f(t)]_p$,
as the \emph{derivative of $F$}. For instance, Fig. 18 shows the derivatives $D(\widetilde{H}^0_1)$ and $D(\widetilde{H}^a_1)$:
\begin{figure}[!ht]
\includegraphics[scale=0.45]{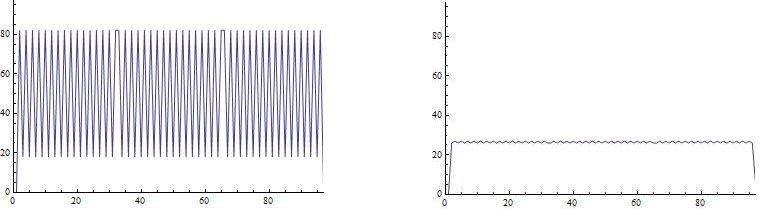}
\caption{$D(\widetilde{H}^0_1(t))$ (left) and $D(\widetilde{H}^a_1(t))$ (right) for $t \in \mathbb{Z}^{\ast}_{97}$ with $a=12$. \label{overflow}}
\end{figure}

\hspace{-0.6cm} Figure 19 shows the graphs of the derivatives $D(\widetilde{H}^0_2)$ and $D(\widetilde{H}^a_2)$:

\begin{figure}[!ht]
\includegraphics[scale=0.45]{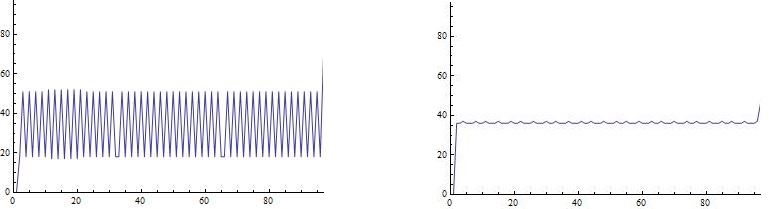}
\caption{$D(\widetilde{H}^0_2(t))$ (left) and $D(\widetilde{H}^a_2(t))$ (right) for $t \in \mathbb{Z}^{\ast}_{97}$ with $a=12$. \label{overflow}}
\end{figure}

\hspace{-0.6cm} What we see is that the functions $\widetilde{H}^0_i, \widetilde{H}^a_i$ are semilinear, that is they have semi-constant first derivative. In particular, this derivative uniformly changes for $i = 1, 2$. It is now understandable to expect, that such a collection of four semi-linear curves in the plane $\mathbb{Z}_p^2$ cannot have a
common solution. Showing this, in general, however, requires some further analysis.


\begin{thebibliography}{10}

\bibitem{FLT} I. Kleiner. \newblock From Fermat to Wiles: Fermat's Last Theorem Becomes a Theorem. \newblock Elem. Math. 55: 19--37, 2010.

\bibitem{W} A. ~Wiles. \newblock Modular elliptic curves and Fermat's Last Theorem. \newblock Annals of Mathematics. 142 (3), 443--551, 1995.

\end{thebibliography}
\end{document}